\documentclass{article}
\usepackage{amsmath}
\usepackage{amsfonts}
\usepackage{amsthm}
\usepackage{centernot}
\usepackage{url}

\numberwithin{equation}{section}

\def\congruent{\equiv}
\def\mod{\bmod}

\def\notmid{\centernot\mid}
\def\ratQ{\mathbb{Q}}

\def\scriptp{ \mathfrak{p} }

\newtheorem{Theorem}{Theorem}[section]
\newtheorem{Observation}{Observation}[section]
\newtheorem{lemma}{Lemma}[section]

\def\Jacobi#1#2{{\left( {#1 \over #2} \right)}}
\def\Legendre#1#2{\left( {#1 \over #2} \right)}

\def\notmid{\not\vert\ }

\begin{document}

\begin{center}
{\large\bf 
Elementary treatment of $p^a \pm p^b + 1 = z^2$ 
}

\bigskip

Reese Scott

\end{center}


 revised 19 May 2023

\bigskip

\begin{abstract}  
 We give a shorter simpler proof of a result of Szalay on the equation $2^a + 2^b + 1 = z^2$.  We give an elementary proof of a result of Luca on the equation of the title for prime $p > 2$.  The elementary treatment is made possible by a lemma which is also useful for other Diophantine equations.  
\end{abstract}

MSC: 11D61, 11D99

\section{Introduction}

All solutions to the title equation for prime $p$ and positive integers $a$, $b$, and $z$ have been found by Szalay \cite{Sz} and Luca \cite{Lu}.  

The purpose of this paper is to give shorter and simpler proofs of these results.  In particular, for the case $p > 2$, we can make the treatment completely elementary, eliminating the use of lower bounds on linear forms in logarithms and the use of results in \cite{BHV}. 

Before proceeding, we give a brief discussion of these changes in the proofs in \cite{Lu} and \cite{Sz}, which deal with the title equation 
where $p$ is a prime and $z$, $a$, and $b$ are positive integers.   
Luca \cite{Lu} handles the case $p>2$ using lower bounds on linear forms in logarithms (see \cite[pp.~7--11]{Lu}) and the well known recent work of Bilu, Hanrot, and Voutier \cite{BHV} (see \cite[pp.~12--14]{Lu}).  In Section 3 we obtain short and elementary proofs of Luca's results, without interfering with the clever use of continued fractions in \cite[equation (18)]{Lu}, by using two elementary lemmas which replace the use of linear forms in logarithms and \cite{BHV} 
(see Lemmas 3.1 and 3.2 in Section 3).  The second of these, Lemma 3.2, has a further application given in Section 4:  we establish a bound on $n$ in the equation 
$x^2 + C = y^n$  when $x$ and $y$ are primes or prime powers.  The bound depends only on the primes dividing $C$.  Beukers \cite{Bk} established a bound on $n$ for more general $x$ when $y=2$, and Bauer and Bennett \cite{BB} greatly improved this bound as well as allowing $y$ to take on many specific values.  The bounds of \cite{Bk} and \cite{BB} depend on the value of $y$ and the specific value of $C$.  See also earlier results of Nagell \cite{N} and Ljunggren \cite{Lj}.  We also apply Lemma 3.2 to a result of Cao \cite{Cao}.   

Further, in proving Theorem 1.4 of this paper, we have removed Luca's use of work of Carmichael \cite{Car}.  Gary Walsh pointed out to the author that \cite{Car} is not needed for proving an auxiliary lemma used by Luca to prove Theorem 1.5 of this paper; although this auxiliary lemma is not used in our proof of Theorem 1.5, Walsh's comment led to our new proof of Theorem 1.4.   
 
Szalay \cite{Sz} handles the equation $2^r - 2^s + 1 = z^2$ using a non-elementary bound of Beukers \cite{Bk}.  However, an earlier result of Beukers, the elementary Theorem 4 of \cite{Bk1}, can be used instead, making Szalay's result elementary, so we will not need to give a new proof in this case.  Szalay \cite{Sz} also handles the case $2^r + 2^s + 1 = z^2$ using a non-elementary result in \cite{Bk}.  In this case we have not obtained a strictly elementary proof; however, we do give a shorter proof of Szalay's result for the case $2^r + 2^s + 1 = z^2$ by replacing the older bound in \cite{Bk} with the recent sharp result of Bauer and Bennett \cite{BB}, not available to Szalay.  Szalay's proof can be further shortened by observing that the methods of his Lemma 8 alone suffice to give the desired contradiction to Beukers' (or Bauer and Bennett's) results; the remaining auxiliary results in \cite{Sz}, including the mapping of one set of solutions onto another, are of independent interest.  An outline of a proof of this result was also given by Mignotte; see the comments at the end of Section D10 of \cite{G}.  

The relevant results of Szalay and Luca are the following:  

\begin{Theorem}[Szalay]   
The equation   
$$  2^a + 2^b + 1 = x^2  \eqno{(1.1)}$$  
has no solutions in positive integers $(a,b,x)$ with $a \ge b $ except for the following cases:  
$$ (a,b,x) = (2t, t+1, 2^t+1) { \rm \  for\ positive\ integer\ } t \eqno{(A)} $$  
$$  (a,b,x) = ( 5,4,7) \eqno{(B)}$$  
$$  (a,b,x) = (9,4,23)  \eqno{(C)}$$  
\end{Theorem}

\begin{Theorem}[Szalay] 
The equation   
$$  2^a - 2^b + 1 = x^2  \eqno{(1.2)}$$  
has no solutions in positive integers $(a,b,x)$ with $a>b $ except for the following cases:  
$$ (a,b,x) = (2t, t+1, 2^t-1) { \rm for\ positive\ integer\ } t>1 \eqno{(D)} $$  
$$  (a,b,x) = ( 5,3,5) \eqno{(E)}$$  
$$  (a,b,x) = (7,3,11)  \eqno{(F)}$$  
$$  (a,b,x) = (15,3,181)  \eqno{(G)}$$
\end{Theorem}

\begin{Theorem}[Luca]  
The only solutions of the equation   
$$p^a \pm p^b + 1 = x^2 \eqno{(1.3)}$$  
in positive integers $(x, p, a, b)$, with $a>b$, and $p$ an odd prime number are  
$(x,p,a,b) = (5,3,3,1)$, $(11,5,3,1)$.   
\end{Theorem} 
  
Luca divides this theorem into three subsidiary theorems:   
  
\begin{Theorem}[Luca]  
The equation   
$$x^2 = y^a + \varepsilon_1 y^b + \varepsilon_2, \qquad \varepsilon_1, \varepsilon_2 \in \{1, -1\}, \eqno{(1.4)}$$  
has no positive integer solutions $(x,y,a,b)$ with $a>b$, $a$ even, and $y>2$ and not a perfect power of some other integer.   
\end{Theorem}    

\begin{Theorem}[Luca]  
There are no solutions to the equation  
$$p^a + p^b + 1 = x^2 \eqno{(1.5)}$$  
in positive integers $(x,p,a,b)$ with $p$ an odd prime.     
\end{Theorem}  

In \cite{Lu} Luca assumes $a$ is odd, but his proof, with only slight changes, gives the more general version used here as Theorem 1.5. 
  
\begin{Theorem}[Luca]  
The only solutions to the equation  
$$p^a - p^b + 1 = x^2 \eqno{(1.6)}$$  
in positive integers $(x,p,a,b)$ with $a>b$ and $p$ an odd prime are   
$(x,p,a,b) = (5,3,3,1)$, $(11,5,3,1)$.  
\end{Theorem}     
  
We would like to thank Michael Bennett for calling our attention to Szalay's paper, Gary Walsh for calling our attention to Luca's paper, and Robert Styer for invaluable suggestions and assistance in preparing this paper.

\section{A shortened proof of Szalay's result}  

\begin{proof}[Proof of Theorem 1.1:]  
Assume (1.1) has a solution that is not one of (A), (B), or (C).  It is an easy elementary result that the only solution to (1.1) with $a=b$ is given by Case (A) with $t = 1$, so we can assume hereafter $a>b$.  

Considering (1.1) modulo 8, we get $b > 2$.  If $b=3$, then $2^a = x^2 - 2^3 - 1 = (x+3) (x-3)$, giving $x = 5$, which is Case (A) with $t=2$, so we can assume hereafter $b>3$.   
  
Write $x = 2^t k \pm 1$ for $k$ odd and the sign chosen to maximize $t>1$.  In what follows, we will always take the upper sign when $x \congruent 1 \bmod{4}$ and the lower sign when $x \congruent 3 \bmod{4}$.      
  
We have   
$$2^a + 2^b + 1 = 2^{2t} k^2 \pm (k \mp 1) 2^{t+1} + 2^{t+1} + 1. \eqno{(2.1)}$$  
From this we see $b=t+1$ so that $t \ge 3$.  Now (2.1) yields $a \ge 2t-1$ with equality only when $t=3$, $k=1$, and $x \congruent 3 \bmod 4$, which is Case (B), already excluded.  So $a \ge 2t$, hence $a > 2t$ since Case (A) has been excluded.  So now   
$$ k \mp 1 = 2^{t-1} g  {\rm \  for\ some\ odd\ } g>0  $$  
We have  
$$  2^{a-2t} = k^2 \pm g = 2^{2t-2} g^2  \pm 2^t g + 1 \pm g  \eqno{(2.2)}$$ 
(2.2) yields $a-2t \ge 2t - 3$ with equality only when $t=3$, $g=1$, and $x \congruent 3 \bmod 4$, which is Case (C), already excluded.  So now $g \pm 1 = 2^t h$ for some odd $h>0$. So we must have $g \ge 2^t \mp 1$.  Assume $x \congruent 3 \bmod{4}$.  Then from (2.2) we derive   
$$  2^{a-2t} > g^2 (2^{2t-2} -1) >  2^{2t} 2^{2t-3} = 2^{4t -3}  \eqno{(2.3)}$$  
Now assume $x \congruent 1 \bmod{4}$.  Then  
$$ 2^{a-2t} > 2^{2t-2} g^2 \ge 2^{2t-2} ( 2^{2t} - 2^{t+1} + 1) > 2^{2t-2} 2^{2t-1} = 2^{4t-3}  $$  
In both cases we have   
$$ a \ge 6t - 2 = 6 b - 8  \eqno{(2.4)}$$    
  
Now we can use Corollary 1.7 in Bauer and Bennett \cite{BB}:   
$$a < {2 \over 2 - 1.48} { \log(2^b+1) \over \log(2) }$$  
Thus,   
$$a < {1 \over 0.26} { \log(2^b+1) \over \log(2^b) } b  < { 1 \over 0.26} {\log(17) \over \log(16) }  b  < 4b$$   
Combining this with (2.4) we obtain $b < 4$ which is impossible since $b>3$.  
\end{proof}      
  
A similar treatment handles Theorem 1.2, although here we must use the familiar results on the equation $x^2 + 7 = 2^y$ to handle the case $b=3$, and also use a slightly more refined computation to establish the second inequality of (2.3), which here applies to $x \congruent 1 \bmod{4}$.  As pointed out in the introduction, however, Szalay already has a short proof of Theorem 1.2 which can be made elementary.

\section{Elementary proofs of Luca's results}   

\begin{proof}[Proof of Theorem 1.4:]  
First we consider the case $b$ even.    
We establish some notation as in \cite{Lu}.  Letting $X=x$, $Y=y^{b/2}$, and $D = y^{a-b}+\varepsilon_1$, we rewrite (1.4) as   
$$X^2 - D Y^2 = \varepsilon_2. \eqno{(3.1)}$$  
The least solution of $U^2 - D V^2 = \pm 1$ is $(U, V) = (y^{(a-b)/2}, 1)$.  Write $X_n + Y_n \sqrt{D} = (y^{(a-b)/2} + \sqrt{D})^n$ for any integer $n$.  For some $j > 1$, $(X, Y) = (X_j, Y_j)$.  As in \cite{Lu}, it is easily seen that $2 | j$.  At this point we diverge from \cite{Lu} and apply Lemmas 1--3 of \cite{Sc} to see that, if $j > 2$, there exists a prime $q$ such that $q | y$, $q | (Y_j/Y_2)$, $Y_{2q} | Y_j$, and $Y_{2q}/ (q Y_2)$ is an integer prime to $y$.  But since $Y_{2q} / (q Y_2)$ is greater than 1 and divides $Y_j$, we have a contradiction.  So $j = 2$ and we must have 
$$y^{b/2} = Y = Y_2 = 2 y^{(a-b)/2}. \eqno{(3.2)}$$ 
  
Now we consider the case $b$ odd and again establish notation as in \cite{Lu}.  Letting $X=x$, $Y = y^{(b-1)/2}$, and $D = y ( y^{a - b} + \varepsilon_1)$, we rewrite (1.4) as (3.1).  At this point we diverge from \cite{Lu} and apply an old theorem of St\"ormer \cite{Sto}: his Theorem 1 says if every prime divisor of $Y$ divides $D$ in (3.1), then $(X, Y) = (X_1, Y_1)$, the least solution of (3.1).  Theorem 1 of \cite{Sto} also applies to show that $(2 y^{a-b} + \varepsilon_1,  2 y^{(a -b - 1)/2} )$ is the least solution $(U_1, V_1)$ of $U^2 - D V^2 = 1$.  If $\varepsilon_2 = -1$, then $2 X_1 Y_1 = 2 y^{(a-b-1)/2}$, which is impossible since $(X_1, y) = 1$, and $y>2$ implies $x = X_1 > 1$.  Thus we must have $\varepsilon_2 = 1$, so that 
$$y^{(b-1)/2} = Y = Y_1 = V_1 = 2 y^{(a -b - 1)/2}. \eqno{(3.3)}$$
At this point we return to \cite{Lu} where it is pointed out that (3.2) and (3.3) require $y=2$ which is not under consideration.
\end{proof}   

We note that Theorem 1 of \cite{Sto} has a short elementary proof.

For Theorems 1.5 and 1.6 we will need the following:

\begin{lemma}   
Let $D$ be any squarefree integer, let $u$ be a positive integer, and let $S$ be the set of all numbers of the form $r + s \sqrt{D}$, where $r$ and $s$ are nonzero rational integers, $(r,sD)=1$, and $u | s$.  Let $p$ be any odd prime number, and let $t$ be the least positive integer such that $\pm p^t $ is expressible as the norm of a number in $S$, if such $t$ exists.  Then, if $\pm p^n$ is also so expressible, we must have $t | n$.  (Note the $\pm$ signs in the statement of this lemma are independent.)   
\end{lemma} 

\begin{proof}   
Assume that for some $p$ and $S$, there exists $t$ as defined in the statement of the lemma.  Then $p$ splits in $\ratQ(\sqrt{D})$; let $[p] = P P'$.  For each positive integer $k$ there exists an $\alpha$ in $S$ such that $P^{kt} = [\alpha]$.  Now suppose $\pm p^{kt+g}$ equals the norm of $\gamma$ in $S$ where $k$ and $g$ are positive integers with $g < t$.  Since $P^{kt+g} $ must be principal, $P^{g} = [\beta]$ for some irrational integer $\beta \in \ratQ(\sqrt{D})$.  Therefore, for some unit $\epsilon$, either $\gamma = \epsilon \alpha \beta$ or $\bar{\gamma} = \epsilon \alpha \beta$.  $\epsilon \alpha \beta$ has integer coefficients and the norm of $\alpha$ is odd, so $\epsilon \beta$ has integer coefficients.  Now $\alpha \in S$ and $\epsilon \alpha \beta \in S$, so that one can see that $\epsilon \beta \in S$, which is impossible by the definitions of $t$ and $g$.  
\end{proof} 

REMARK: If $u=1$, the proof of Lemma 3.1 can be shortened and simplified by noting that $t$ must be the least integer such that $P^t$ is principal (so that $P^g$ cannot be principal, immediately giving $t \mid n$), citing Observation 3.1 which follows.  

\begin{Observation} 
Let $D$ be a positive squarefree integer.  Any principal ideal $A$ in $\ratQ(\sqrt{D})$ having an odd norm has a generator with rational integer coefficients.  
\end{Observation}

\begin{proof}
The Observation clearly holds for $D \not\equiv 1 \bmod 4$.  If $D \equiv 1 \bmod 8$, then any integer of the field having an odd norm must have rational integer coefficients.  So we can assume $D \equiv 5 \bmod 8$.  Now Observation 3.1 follows from Observation 3.2 which follows.
\end{proof}

\begin{Observation} 
Let $D \equiv 5 \bmod 8$ be a positive squarefree integer.  Then, for any integer $\alpha \in \ratQ(\sqrt{D})$ there exists a unit $\delta \in \ratQ(\sqrt{D})$ such that $\alpha \delta$ has rational integer coefficients. 
\end{Observation}

\begin{proof}  Observation 3.2 follows from combining (i) and (ii) below.  

(i) Let $\beta = \frac{r + s \sqrt{D}}{2}$ where $2 \notmid rs$ and $ r \not\equiv s \bmod 4$.  Let $\gamma = \frac{h+k \sqrt{D}}{2}$ where $2 \notmid hk$ and $h \equiv k \bmod 4$.  
$$\beta \gamma = \frac{ (rh + skD) + (rk + sh) \sqrt{D} }{4}. $$
$rh + skD \equiv rh + sk \equiv h(r+s) \equiv 0 \bmod 4$ and $rk+sh \equiv h(r+s) \equiv 0 \bmod 4$, so that $\beta \gamma$ has rational integer coefficients.  

(ii) Let $\delta$ be any integer in $\ratQ(\sqrt{D})$ with $\delta = \frac{u + v \sqrt{D}}{2}$ where $2 \notmid uv$.  
Then $\delta^2 = \left( \frac{u+v \sqrt{D}}{2} \right)^2 = \frac{u_2 + v_2 \sqrt{D}}{2}$ where 
$2 \notmid u_2 v_2$ and $u$, $v$, $u_2$, $v_2$ have the following property: if $u$ is congruent (respectively, not congruent) to $v$ modulo 4, then $u_2$ is not congruent (respectively, congruent) to $v_2$ modulo 4.  To see this, note that
$$ u_2 = \frac{ u^2 + v^2 D }{2}, v_2 = uv,$$
so that $u_2 \equiv 3 \mod 4$ (recall $D \equiv 5 \bmod 8$).  
\end{proof} 

This concludes the handling of the Remark following Lemma 3.1 above.

\begin{proof}[Proof of Theorem 1.5:]    
We first establish some notation by paraphrasing \cite[Section 3]{Lu}:   
Looking at (1.5), we see that the only case in which solutions might exist is when $p \congruent 3 \bmod 4$ and $2 \notmid a-b$; choose $b$ even (without loss of generality).   Let $p^b + 1 = D u^2$, with $D$ square-free and $u > 0$ an integer.  At this point we diverge from \cite{Lu} and note that if $S$ is the set of all integers of the form $r + s \sqrt{D}$ with nonzero rational integers $r$ and $s$, $(r,sD) = 1$ and $u | s$, then $p^a$ and $-p^b$ are both expressible as the norms of numbers in $S$.  Therefore Lemma 3.1 shows that $\pm p^c$ is expressible as the norm of a number in $S$, where $c$ divides both $a$ and $b$.  From this  point on, we return to the method of proof of \cite{Lu}: $a$ is odd and $b$ is even, so we have $c \le b/2$.  For some coprime positive integers $v$ and $w$ such that $(v, p^b+1) = 1$, we must have  
$$  v^2 - w^2 (p^b+1) = \pm p^c. \eqno{(3.4)}$$  
(3.4) corresponds to (17) in \cite{Lu}.  Since $| p^c| < \sqrt{p^b+1}$, $v/w$ must be a convergent of the continued fraction for $\sqrt{p^b+1}$.  But then, since $p^b + 1$ is of the form $m^2 + 1$, we must have $p^c = \pm 1$, impossible.
\end{proof}

\begin{proof}[Proof of Theorem 1.6:]     
As in \cite{Lu}, we write $p^b - 1 = D u^2$, $D$ and $u$ positive integers and $D$ squarefree, and consider the equation   
$$p^n = h^2 + k^2 u^2 D \eqno{(3.5)}$$ 
in relatively prime nonzero integers $h$ and $k$, and positive integer $n$.   
From (1.6) we see that (3.5) has the solutions $(n, h, k) = (b, 1,1)$ and $(a, x, 1)$.  Clearly, $p$ splits in $\ratQ(\sqrt{-D})$, and we can let $[p] = \pi_1 \pi_2$ be its factorization into ideals.  We can take  
$${\pi_1}^b  = [ 1 + u \sqrt{-D} ], {\pi_1}^a = [x \pm u \sqrt{-D} ].  \eqno{(3.6)}$$
At this point we diverge from \cite{Lu}: clearly $b$ is the least possible value of $n$ in (3.5), so we can apply Lemma 3.1 to obtain $b | a$.  Thus,   
$$( 1 + u \sqrt{-D} )^{a/b} = ( x \pm u \sqrt{-D} ) \epsilon \eqno{(3.7)}$$  
where $\epsilon$ is a unit in $\ratQ(\sqrt{-D})$.  If $D = 1$ or 3, we note $2 | u$ and $ 2 \notmid x$, so that we must have $\epsilon = \pm 1$.   Thus, using (3.7), we see that Theorem 1.6 follows immediately upon establishing the following elementary lemma (Lemma 3.2).  
\end{proof}  
  
\begin{lemma} 
The equation   
$$(1+ \sqrt{-D} )^r = a \pm \sqrt{-D} \eqno{(3.8)}$$  
has no solutions with $r>1$ when $D$ is a positive integer congruent to 2 mod 4 and $a$ is any integer, except for $D=2$, $r=3$.    
  
Further, when $D$ congruent to 0 mod 4 is a positive integer such that $1+D$ is prime or a prime power, (3.8) has no solutions with $r>1$ except for $D=4$, $r=3$.   
\end{lemma}  
  
(Note that here $D$ corresponds to $D u^2$ in the proof of Theorem 1.6 above, which follows the notation of \cite{Lu}.  Thus, in the proof of Lemma 3.2, $D$ is not necessarily squarefree.  Note also $r$ corresponds to $a/b$ in the proof of Theorem 1.6.)    
  
\begin{proof}[Proof of Lemma 3.2:]  
Assume (3.8) has a solution with $r>1$ for some $a$ and $D$.  From Theorem 13 of \cite{BH}, we see that, if $r>1$, then $r$ is a prime congruent to 3 mod 4 and there is at most one such $r$ for a given $D$.  Thus we obtain   
$$(-1)^{{D+2 \over 2}} = r - {r \choose 3} D + {r \choose 5} D^2 - \dots - D^{{r-1 \over 2}} \eqno{(3.9)}$$  
If $r=3$, (3.9) shows that $|D-3| = 1$, giving the two exceptional cases of the Lemma.  So from here on we assume $3 \notmid r$.    
  
We will use two congruences:   
$$  (-1)^{{D+2 \over 2}} \congruent \Legendre{ r}{ 3} 2^{r-1} \bmod D-3 \leqno{{\rm Congruence\ 1:}}$$  
$$  (-1)^{{D+2 \over 2}} \congruent 2^{r-1} \bmod D+1  \leqno{{\rm Congruence\ 2:}}$$  
Congruences 1 and 2 correspond to congruences (9e) and (9f) of Lemma 7 of \cite{BH}.  From Congruence 1 we see that $D-3$ cannot be divisible both by a prime 3 mod 4 and a prime 5 mod 8.  So $D \congruent 2 \bmod 4$ implies $D \not\congruent 3 \bmod 5$.  Now let $D+1 = y$.  If $D \congruent 1 \bmod 5$, $y^r \congruent 3 \bmod 5$; since $a^2 + D = y^r$, $a^2  \congruent 2 \bmod 5$, impossible.  If $D \congruent 2 \bmod 5$, $y^r \congruent 2 \bmod 5$, so that 5 divides $a$.  Since in this case  $D$ is a quadratic nonresidue modulo 5, we see from (3.8) that  $5 | a$ implies $3 | r$, which we have excluded.  Now $y^r$ is congruent to $-y$ modulo   
$y^2 + 1$ so that $a^2$ is congruent $-2 y + 1$ modulo $y^2 + 1$.  So, using the Jacobi symbol, we must have   
$$ 1 = \Jacobi{-2y+1}{(y^2+1)/2} = \Jacobi{2y^2+2}{ 2y-1} = \Jacobi{y+2}{ 2y-1} = \Jacobi{-5}{ y+2}$$  
If $D \congruent 2 \bmod{4}$, then $y \congruent 3 \bmod 4$ and the last Jacobi symbol in this sequence equals $\Jacobi{y+2}{ 5} = \Jacobi{D+3}{5}$, which has the value -1 when $D$ is congruent to 0 or 4 modulo 5.  Thus, when $D \congruent 2 \bmod 4$, we have shown that there are no values of $D$ modulo 5 that are possible.    
  
So we assume hereafter that $D \congruent 0 \bmod 4$.  Write $D+1 = p^n$ where $p$ is prime, and let $g$ be the least number such that $2^g \congruent -1 \bmod p$, noting Congruence 2.  We see that $g | r-1$ and also $g | p-1 | p^n - 1 = D$.  Now (3.9) gives $-1 \congruent 1 \bmod g$ so that $g \le 2$.  Assume first that $n$ is odd.  Since $4 | D$, $p \congruent 1 \bmod 4$.  In this case, we must have $g=2$, $p=5$.  If $n$ is even, since we have $1+D = p^n$ and $a^2 + D = p^{rn}$,  we must have $2 p^{rn / 2} - 1 \le D = p^n - 1$, giving $r < 2$, impossible.  So we have $n$ odd, $p=5$.    
  
Since $n$ is odd, $D \congruent 4 \bmod 8$, and, since $r \choose 3$ is odd, (3.9) gives $r \congruent 3 \bmod 8$.  Now assume $r \congruent 2 \bmod 3$ and let $y = 5^n = 1+D$.  Then $y^r \congruent y^2 \bmod y^3 - 1$, so that $a^2 \congruent y^2 -y + 1 \bmod y^2 + y + 1$, so that   
$$1 = \Jacobi{y^2 - y + 1}{ y^2 + y + 1} = \Jacobi{-2y}{ y^2 +y+1} = \Jacobi{-2}{ y^2 + y + 1}$$  
which is false since $y^2 + y+1 \congruent 7 \bmod 8$.  Thus we have $r \congruent 19 \bmod 24$ so that $y^r \congruent -y^7 \bmod y^{12} + 1$, so that $a^2 \congruent -y^7 - y + 1 \bmod {y^{12} + 1 \over 2}$.  Thus we have   
\begin{align*} 
 1 &= \Jacobi{-y^7 - y + 1}{ (y^{12}+1)/2} = \Jacobi{y^7 + y - 1}{ (y^{12} + 1)/2} = \Jacobi{ 2 (y^{12} + 1)}{ y^7 + y- 1} \\
&= \Jacobi{y^{12} + 1}{ y^7 + y -1} = \Jacobi{ y^6 - y^5 - 1}{ y^7 + y - 1} 
= \Jacobi{ y^7 + y -1}{ y^6 - y^5 - 1} \\
&= \Jacobi{ y^5 + 2 y }{ y^6 - y^5 - 1 } = \Jacobi{y^4 + 2}{ y^6 - y^5 - 1}  
= - \Jacobi{y^6 - y^5 - 1}{ y^4 + 2} \\ 
&=  \Jacobi{2 y^2 -2y+1}{ y^4 + 2} = \Jacobi{ y^4 + 2}{ 2 y^2 - 2y + 1} = \Jacobi{7}{ 2 y^2 - 2y + 1} \\
&= \Jacobi{ 2 y^2 - 2y + 1}{ 7}    \\
\end{align*}
which is possible only when $y$ is congruent to 1, 4, or 0 modulo 7.  This is impossible since   
$y$ is an odd power of 5.  
\end{proof}

\section{Further applications of Lemma 3.2 }  

In this section we show how Lemma 3.2 can be used to handle Diophantine equations other than (1.6). We give two examples of such equations.  The first of these (Equation (4.1) below) is the familiar Ramanujan-Nagell equation for the special case $x$ and $y$ primes or prime powers; we obtain a bound on $n$ in (4.1) where the bound depends only on the primes dividing $C$.   The second of these (Equation (4.6) below) appears in Theorem 4.2 below, which is a familiar result of Cao; we use Lemma 3.2 to simplify a proof of Cao's result given in \cite{ScSt} (as far as we know, Cao did not publish a proof, although existence of a proof is mentioned in \cite{Cao} which is an abstract of a paper in Chinese).

\begin{Theorem} 
\label{Thm4}
Let $C$ be a positive integer, and let $PQ$ be the largest squarefree divisor of $C$, where $P$ is chosen so that $(C/P)^{1/2}$ is an integer.  If the equation   
\begin{equation}  x^2 + C = y^n \label{(4.1)}\end{equation}  
has a solution $(x,y,n)$ with $y$ a prime and $x$ divisible by at most one prime, $(x,y)=1$, $n$ a positive integer, and $(x,y,n) \ne (401,11,5)$, $(7,3,4)$, $(1,2,2)$, $(5,2,5)$, $(11,2,7)$, $(181,2,15)$, $(11,2,8)$, $(61,2,12)$, or $(13, 2, 9)$, 
then we must have either $n=3$ or   
$$(n-s) \mid  N = 2 \cdot 3^\ell h(-P) \langle q_1 - \Legendre{-P}{q_1}, \dots, q_n - \Legendre{-P}{q_n} \rangle$$  
where $s=2$ or $0$ according as $y=2$ or not, $\ell=1$ or $0$ according as $3 < P \congruent 3 \bmod 8$ or not, $h(-P)$ is the lowest $m$ such that $\mathfrak{a}^m$ is principal for every ideal $\mathfrak{a}$ in $\ratQ(\sqrt{-P})$, $\langle a_1, a_2, \dots, a_n \rangle$ is the least common multiple of the members of the set $S = \{ a_1, a_2, \dots, a_n \}$ when $S \ne \emptyset$, $\langle a_1, a_2, \dots, a_n \rangle = 1$ when $S = \emptyset$, $q_1 q_2 \dots q_n = Q$ is the prime factorization of $Q$, and $\Legendre{a}{q}$ is the familiar Legendre symbol unless $q=2$ in which case $\Legendre{a}{2} = 0$. 
\end{Theorem} 
  
\begin{proof}
Assume there exists a solution to (\ref{(4.1)}).  Let $\scriptp \bar{\scriptp}$ be the prime ideal factorization of $y$ in $\ratQ(\sqrt{-P})$.  

\bigskip
\noindent
{\it Case 1}:  $2 \mid C$.  Let $k$ be the smallest number such that $\scriptp^k = [\alpha]$ is principal with a generator $\alpha$ having integer coefficients.  When $P=1$, we choose $\alpha$ so that the coefficient of its imaginary term is even.  When $P = 3$ we can take $k=1$. Then  
$$\alpha^{n/k} = \pm x \pm \sqrt{-C}$$  
where the $\pm$ signs are independent.  Note that when $P=3$ and $\alpha^{n/k} \epsilon = x \pm \sqrt{-C} $ for some unit $\epsilon$, we must have $\epsilon = \pm 1$.  Let $j$ be the least number such that $\alpha^j = u + v Q \sqrt{-P}$ for some integers $u$ and $v$.  By elementary properties of the coefficients of powers of integers in a quadratic field, $jk  | N/2$.  Also, $jk | n = jkr$ for some $r$.    
So we have  
$$(u + v Q \sqrt{-P})^r = \pm x \pm \sqrt{-C} $$  
If $r=1$ or $r=2$, Theorem~\ref{Thm4} holds, so assume $r \ge 3$.    
  
If $r$ is even, then any prime dividing $u$ must divide $C$, since $\pm x \pm \sqrt{-C}$ must be divisible by $( u  + v Q \sqrt{-P})^2.$  Since $(u, C) = 1$, we must have $u = \pm 1$ when $r$ is even.    
  
If $r$ is odd, then $u$ divides $x$.  $x = \pm 1$ implies $u = \pm 1$.  Assume $| x | > 1$.  Let $x = \pm g^b$ where $g$ is a positive prime and $b > 0$.  Then, when $r$ is odd, $u = \pm g^t$ for some $t \ge 0$  (this follows from the same kind of elementary reasoning used for Lemmas 1--3 of \cite{Sc}).  Also, every prime dividing $v$ divides $C$.  Thus, if $t>0$, then by Theorem 1 of \cite{Sc}, $r=1$ which we already excluded.   (Note that the only relevant exceptional case in Theorem 1 of \cite{Sc} is $(x,y,C) = (3,13,10)$, in which case $n = 1$ or $3$.)   
  
So $u = \pm 1$ regardless of the value of $x$ or the parity of $r$.  Letting $D = v^2 Q^2 P$, we have  
$$(1+ \sqrt{-D})^r = \pm x \pm w \sqrt{-D} $$  
for some positive integer $w$.  If $w=1$, we see from Lemma~3.2 that $r=3$ and $j = k = 1$, so that $n = 3$ and the theorem holds.    
  
So $w>1$, and $w$ is divisible only by primes dividing $C$.  In what follows, we apply Lemmas 1--3 of \cite{Sc}.  We must have at least one prime $r_1$ dividing $C$ which also divides $r$.  We have, for any such $r_1$,  
\begin{equation}(1+\sqrt{-D} )^{r_1} = \pm x_1 \pm w_1 \sqrt{-D} \label{(4.2)}\end{equation}  
where $w_1 | w$.  If $r_1$ is odd, we have  
\begin{equation} \pm w_1= r_1 - {r_1 \choose 3} D + {r_1 \choose 5} D^2 - \dots + (-D)^{(r_1 -1)/ 2}. \label{(4.3)}\end{equation}   
From (\ref{(4.3)}) we see that $r_1 \mid w_1$ and $w_1$ is divisible by no primes other than $r_1$.  If $r_1 > 3$, then $r_1^2 \notmid w_1$, so that $w_1 = \pm r_1$.    
  
If $r_1 = 3$, we must have $w_1 = \pm 3^z$ for some $z > 0$ so that $D= 3^z + 3$.  Now $1+D$ is the norm of $\alpha^{j}$ which equals $y^{jk}$.  But $1+D = 3^z + 4$ cannot be a perfect power of $y$ by Lemma 2 of \cite{ScSt}.  So $j=k=1$.  Now $| x_1 | = 3D-1 > 1$.  Also, $(x_1, C) = 1$ so $2  \notmid \frac{r}{r_1}$.  Thus, $|x_1|$ must be a power of the prime dividing $x$ (this follows from the same kind of elementary reasoning used for Lemmas 1--3 of \cite{Sc}).  By Theorem 1 of \cite{Sc}, $r= r_1$, $n = 3jk = 3$, and the theorem holds.    
  
If $r_1 = 5$ then (\ref{(4.3)}) shows that $\pm 5 = 5 - 10 D + D^2$.  Since $5 | D$, this implies $D=10$, $y^{jk}=11$ which gives $(x_1, y, r_1, j,k) = (401, 11, 5, 1, 1)$.  If $r> r_1$, we must have $2 \notmid r$ and $401 | x$, so Theorem 1 of \cite{Sc} shows $r=r_1$. This leads to the case $(x,y,n) = (401, 11, 5)$.    
  
If $r_1 \ge 7$, (\ref{(4.3)}) is impossible for $w_1 = \pm r_1$.   
  
Finally, it remains to consider $r = 2^h$, $h > 1$.  Then we have (\ref{(4.2)}) with $r_1 = 2$, $| x_1 | = D-1$.  If $D > 2$, then, since $D-1>1$, we have $2 \notmid \frac{r}{r_1}$, contradicting $h>1$.     
So $D=2$, so that $y^{jk} = 1+ D = 3$, and $n = r = 2^h$.   $n = 4$ gives the exceptional case $(x,y,n) = (7,3,4)$; and $n>4$ gives $7 \mid w$, impossible. 

\bigskip 
\noindent 
{\it Case 2}: $y = 2$.  
The theorem holds for the unique trivial case in which $n=1$, and we can eliminate from consideration the exceptional case $(x,y,n)=(1,2,2)$, so that we can assume $n>2$, so that $P \equiv 7 \bmod 8$.  Let $k$ be the smallest number such that $\scriptp^{k}$ is principal, and let $j$ be the least number such that $\scriptp^{jk} = \left[\frac{u + v Q \sqrt{-P} }{2} \right]$ where $u$ and $v$ are integers.  As in Case 1 we obtain $j k \mid N/2$, $jk \mid (n-2) = jkr$ for some $r \ge 3$, and $u = \pm 1$.  Letting $D = v^2 Q^2 P$, we have 
$$ (1+ \sqrt{-D})^r = (\pm x \pm w \sqrt{-D}) 2^{r-1}$$
for some positive integer $w$.  If $w=1$, then, using Theorem 4 of \cite{Bk1}, we obtain the fourth, fifth, and sixth exceptional cases listed in the formulation of the theorem.  So we can take $w>1$, noting $w$ is divisible only by primes dividing $C$.  

Still proceeding as in Case 1, we can choose a prime $r_1 \mid C$ such that $r_1 \mid r$.  Since $r_1$ is odd we have 
\begin{equation}  \pm 2^{r_1 - 1} w_1 = r_1 - { r_1 \choose 3} D + { r_1 \choose 5} D^2 - \dots + (-D)^{(r_1-1)/2}. \label{(125)} \end{equation}    
(\ref{(125)}) corresponds to (\ref{(4.3)})  in Case 1.  

If $r_1 = 3$ we must have $w_1 = \pm 3^z$ for some $z>0$ so that $D = 4 \cdot 3^z + 3$, so that $3^z + 1= 2^{jk}$, so that $jk = 2$, $z=1$, $D=15$.  As in Case 1, we can use Theorem 1 of \cite{Sc} to get $r=r_1$, and we obtain the seventh exceptional case in the formulation of the theorem (note that the exceptional cases of Theorem 1 of \cite{Sc} do not apply here).   So the theorem holds for Case 2 when $r_1 = 3$.    

If $r_1 = 5$ then (\ref{(125)}) shows that $\pm 80 = 5 - 10 D + D^2$, from which we obtain $D=15$.  As above, we get $r=r_1$, giving the eighth exceptional case in the formulation of the theorem.  So the theorem holds for Case 2 when $r_1 = 5$.

If $r_1 \ge 7$, then $ |w_1| = r_1$, and, using (\ref{(125)}) and the fact that $D = 2^{jk+2} -1$, we get 
\begin{equation}  \frac{2^{r_1-1} -1}{2^{jk+2} -1} = - \frac{(r_1-1)(r_1-2)}{6} + \frac{D}{r_1} \left( {r_1 \choose 5} - \dots \pm D^{(r_1-5)/2} \right)  \label{(126)} \end{equation} 
noting that the term on the left of (\ref{(125)}) must be positive by consideration modulo $r_1^2$.  Since $t = jk+2$ is the least $t$ satisfying $2^t \equiv 1 \bmod D$, we must have $jk+2 \mid r_1 - 1$ and the left side of (\ref{(126)}) is congruent $(r_1-1)/(jk+2)$ modulo $r_1$.  Now considering  (\ref{(126)}) modulo $r_1$ we find $jk \congruent 1 \bmod r_1$.  Since $jk+2 \le r_1 -1$, we have $jk=1$, so that $D = r_1 = r = 7$, and we obtain the final exceptional case listed in the theorem.

\bigskip 
\noindent 
{\it Case 3}:  $2 \mid x$.  
If $C=1$ it is an easy elementary result that we must have $n =1$, and the theorem holds.  If $C>1$, then the proof of the theorem for Case 3 is immediate after Theorem 2 of \cite{Sc}.     
\end{proof}

\begin{Theorem} 
The equation 
$$q^n+2^h = p^w \eqno{(4.6)}$$
where $p$ and $q$ are distinct positive primes, $n$ is even, $h$ is a positive integer, and $w>1$, has as its only solutions 
$3^2+2^4 = 5^2$, $7^2+2^5=3^4$, $5^2+2 = 3^3$, and $11^2+2^2=5^3$. 
\end{Theorem} 

\begin{proof} 
A proof is given in \cite[Lemma 4]{ScSt} which is simplified by using Lemma 3.2 directly instead of using the result of Theorem 13 of \cite{BH} on which Lemma 3.2 depends.
\end{proof}

\end{document}